\documentclass[12pt,reqno]{amsart}
\usepackage{enumerate, latexsym, amsmath, amsfonts, amssymb, amsthm, color}
\def\pmod #1{\ ({\rm{mod}}\ #1)}
\def\Z{\Bbb Z}

\def\l{\left}
\def\r{\right}
\def\bg{\bigg}
\def\({\bg(}
\def\){\bg)}
\def\t{\text}
\def\f{\frac}

\def\exp{{\rm exp}}

\def\ls{\leqslant}
\def\gs{\geqslant}
\def\se {\subseteq}
\def\sm{\setminus}
\def\bi{\binom}
\def\al{\alpha}

\def\eq{\equiv}

\def\Ack{\medskip\noindent {\bf Acknowledgment}}
\theoremstyle{plain}
\newtheorem{theorem}{Theorem}

\newtheorem{lemma}{Lemma}[theorem]
\newtheorem{corollary}{Corollary}[theorem]
\newtheorem{conjecture}{Conjecture}[theorem]
\theoremstyle{definition}

\theoremstyle{remark}

 \vspace{4mm}

\begin{document}
 \baselineskip=17pt
\hbox{Electron. J. Combin. 24(2017), no.\,1, P1.17, 1--6.}
\medskip

\title
[On a permutation problem for finite abelian groups]
{On a permutation problem for finite abelian groups}

\author
[Fan Ge and Zhi-Wei Sun] {Fan Ge and Zhi-Wei Sun*}

\thanks{*This corresponding author is supported by the National Natural Science Foundation of China (grant 11571162)}

 \address {(Fan Ge) Department of Mathematics, University of Rochester
\\Rochester, NY-14627, USA}
\email{{\tt fange.math@gmail.com}}

\address {(Zhi-Wei Sun) Department of Mathematics, Nanjing
University, Nanjing 210093, People's Republic of China}
\email{{\tt zwsun@nju.edu.cn}}


 \begin{abstract} Let $G$ be a finite additive abelian group with exponent $n>1$, and let $a_1,\ldots,a_{n-1}\in G$.
We show that there is a permutation $\sigma\in S_{n-1}$ such that all the elements $sa_{\sigma(s)}\ (s=1,\ldots,n-1)$
are nonzero if and only if
$$\left|\left\{1\ls s<n:\ \f{n}da_s\not=0\right\}\right|\gs d-1\ \ \mbox{for every positive divisor}\ d\ \mbox{of}\ n.$$
When $G$ is the cyclic group $\mathbb Z/n\mathbb Z$, this confirms a conjecture of Z.-W. Sun.
\end{abstract}

\maketitle

\section{Introduction}
\setcounter{theorem}{1}

Let $n\in\Z^+=\{1,2,3,\ldots\}$ and let $S_n$ denote the symmetry group of all permutations on $\{1,\ldots,n\}$.
A conjecture of G. Cramer stated that for any integers $m_1,\ldots,m_n$ with $\sum_{s=1}^nm_s\eq0\pmod n$ there is a permutation $\sigma\in S_n$
such that $1+m_{\sigma(1)},\ldots,n+m_{\sigma(n)}$ are pairwise distinct modulo $n$. In 1952 M. Hall [H] proved an extension of this conjecture.

In 1999 H. S. Snevily [Sn] conjectured that if $n>1$ is an integer and $m_1,\ldots,m_k$ are integers with $k\ls n-1$ then
there is a permutation $\sigma\in S_k$ such that $1+m_{\sigma(1)},\ldots,k+m_{\sigma(k)}$ are pairwise distinct modulo $n$.
This was confirmed by A. E. K\'ezdy and Snevily [KS] in the case $k\ls(n+1)/2$, and an application to tree embeddings was also given in [KS].

Let $n>1$ and $m_1,\ldots,m_{n-1}$ be integers. When is there a permutation $\sigma\in S_{n-1}$ such that none of
the $n-1$ numbers $sm_{\sigma(s)}\ (s=1,\ldots,n-1)$ is congruent to $0$ modulo $n$?
If there is such a
permutation $\sigma$, then for each positive divisor $d$ of $n$ we have
$$|\{1\ls c<d:\ d\nmid m_{\sigma(cn/d)}\}|
\gs\l|\l\{1\ls c<d:\ n\nmid\f {cn}dm_{\sigma(cn/d)}\r\}\r|=d-1,$$
and hence the sequence $\{m_s\}_{s=1}^{n-1}$ has the following property:
\begin{equation}\label{1.1}
\big|\{1\ls s <n:\ d\nmid m_s\}\big|
\gs d-1\ \ \mbox{for any}\ d\in D(n),
\end{equation}
where $D(n)$ denotes the set of all positive divisors of $n$.

In 2004 the second author (cf. \cite{S09}) made the following conjecture.
\begin{conjecture}\label{conj1.1} {\rm (Z.-W. Sun)} Let $n>1$ be an integer. If $m_1,m_2,\ldots,m_{n-1}$ are integers satisfying $(\ref{1.1})$, then
there exists a permutation $\sigma$ on $\{1,\ldots,n-1\}$ such that
$n\nmid sm_{\sigma(s)}$ for all $s=1,\ldots,n-1$.
\end{conjecture}

In this paper we aim to prove an extension of this conjecture for finite abelian groups.

For a finite multiplicative group $G$, its exponent $\exp(G)$ is defined to be the least positive integer such that $x^n=e$ for all $x\in G$, where $e$ is the identity of $G$.
For a finite abelian group $G$, $\exp(G)$ is known to be $\max\{o(x):\ x\in G\}$, where $o(x)$ denotes the order of $x$.
If $G$ is an additive group, then for $k\in\Z^+$ and $a\in G$ we write $ka$ for the sum $a_1+\ldots+a_k$ with $a_1=\cdots=a_k=a$.

\begin{theorem}\label{thm1.1} Let $G$ be a finite additive abelian group with exponent $n>1$.
For any $a_1,\ldots,a_{n-1}\in G$, there is a permutation $\sigma\in S_{n-1}$ such that all the elements $sa_{\sigma(s)}\ (s=1,\ldots,n-1)$
are nonzero if and only if
\begin{align}\label{1.2}
\l|\l\{1\ls s<n:\ \f{n}da_s\not=0\r\}\r|\gs d-1\ \ \mbox{for all}\ d\in D(n).
\end{align}
\end{theorem}

Applying Theorem 2 to the cyclic group $\Z/n\Z$, we immediately confirm Conjecture \ref{conj1.1} of Sun. As an application, we obtain the following result.

\begin{theorem}\label{Thm1.2} Let $m_1,m_2,\ldots,m_{n-1}\ (n>1)$ be integers satisfying $(\ref{1.1})$. Then the set
$$\l\{\sum_{i\in I}m_i:\ I\se\{1,\ldots,n-1\}\r\}$$
contains a complete system of residues modulo $n$.
\end{theorem}

Obviously Theorem 3 extends the following result of the second author (cf. the paragraph following \cite[Theorem 2.5]{S03}).

\addtocounter{theorem}{1}
\begin{corollary}\label{cor1.1} Let $n>1$ be an integer and let $m_1,m_2,\ldots,m_{n-1}$ be integers all relatively prime to $n$.
Then the set $\l\{\sum_{i\in I}m_i:\ I\se\{1,\ldots,n-1\}\r\}$
contains a complete system of residues modulo $n$.
\end{corollary}

As usual, for any $a\in\Z$ and $n\in\Z^+$, we write $(a,n)$ for the greatest common divisor of $a$ and $n$.

Let $n>1$ be an integer. If $m_s\in\Z$ and $(m_s,n)\ls s$ for all $s=1,\ldots,n-1$, then for any $d\in D(n)$ we have
$$|\{1\ls s<n:\ d\nmid m_s\}|\gs |\{1\ls s<n:\ s<d\}|=d-1,$$
and hence by Theorem \ref{thm1.1} for some $\sigma\in S_{n-1}$ we have
$n\nmid \sigma(s)m_s$ for all $s=1,\ldots,n-1$.
This is equivalent to the following theorem in the case $a_1=\cdots=a_{n-1}$.

\begin{theorem}\label{Thm1.3}
Let $m_1,m_2,\ldots,m_{n-1}\ (n>1)$ be integers with $(m_s,n)\ls s$ for all $s=1,\ldots,n-1$.
For any $a_1,\ldots,a_{n-1}\in\Z$,
there is a function $f:\{1,\ldots,n-1\}\to\{1,\ldots,n-1\}$ such that the sums
$$f(1)+a_1,\ \ldots,\ f(n-1)+a_{n-1}$$ are pairwise distinct modulo $n$ and also
none of the numbers $$f(1)m_1,\ \ldots,\ f(n-1)m_{n-1}$$ is divisible by $n$.
\end{theorem}

\addtocounter{theorem}{1}
Motivated by Theorems 2 and 3, we pose the following conjecture.
\begin{conjecture}\label{conj1.2}
Let $G$ be a finite abelian group with exponent $n>1$. If $a_1,\ldots,a_{n-1}$ are elements of $G$
with $sa_s\not=0$ for all $s=1,\ldots,n-1$, then we have
\begin{equation}\bg|\bg\{\sum_{i\in I}a_i:\ I\se\{1,\ldots,n-1\}\bg\}\bg|\gs n.
\end{equation}
\end{conjecture}

By Theorems 2 and 3, this conjecture holds for finite cyclic groups.
For any finite abelian group $G$ with exponent $n>1$, it has a cyclic subgroup $H$ of order $n$,
and hence for $a_1,\ldots,a_{n-1}\in H$ the set $\{\sum_{i\in I}a_i:\ I\se\{1,\ldots,n-1\}\}$
contains at most $n$ elements of $G$.
\medskip

We will show Theorem 2 in the next section and prove Theorems 3 and 5 in Section 3.

\section{Proof of Theorem \ref{thm1.1}}

\medskip

\medskip
\noindent{\it Necessity}: If there is a permutation $\sigma\in S_{n-1}$ such that $sa_{\sigma(s)}\not=0$ for all $s=1,\ldots,n-1$, then
for any $d\in D(n)$ we have
$$\l|\l\{1\ls s<n:\ \f nda_s\not=0\r\}\r|\gs \l|\l\{1\ls c<d:\ \f {cn}da_{\sigma(cn/d)}\not=0\r\}\r|=d-1.$$
This concludes the proof of the necessity. \qed

\medskip
\noindent{\it Sufficiency}:
Suppose, to the contrary, that there are $a_1,\ldots,a_{n-1}\in G$ satisfying (\ref{1.2})
such that the set
$$I(\sigma):=\{1\ls i<n: ia_{\sigma(i)}=0\}=\{1\ls i<n: o(a_{\sigma(i)})\mid i\}$$
is nonempty for any $\sigma\in S_{n-1}$.
Take such $a_1,\ldots,a_{n-1}\in G$ with $\sum_{s=1}^{n-1}o(a_s)$ maximum, and
choose $\sigma\in S_{n-1}$ with $|I(\sigma)|$ minimum.
\medskip

{\it Claim} 1: $|I(\sigma)|=1$.

As $n=\exp(G)$, there is an element $x$ of $G$ with $o(x)=n$.
Let $j\in I(\sigma)$, and
for $s=1,\ldots,n-1$ define
$$a_s^*=\begin{cases} x&\t{if}\ s=\sigma(j),\\a_s&\t{otherwise}.\end{cases}$$
If $(n/d)a_{\sigma(j)}\not=0$ with $d\in D(n)$, then $d>1$ and $(n/d)x\not=0$.
As $o(a_{\sigma(j)})\mid j$, we have $o(a_{\sigma(j)})\ls j<n=o(x)$.
Since $\sum_{s=1}^{n-1}o(a_s^*)>\sum_{s=1}^{n-1}o(a_s)$, by our choice of $a_1,\ldots,a_{n-1}$, for some $\tau\in S_{n-1}$ we have $sa_{\tau(s)}^*\ne 0$ for all $s=1,\ldots,n-1$.
For any $1\ls s<n$ with $\tau(s)\not=\sigma(j)$, we have $sa_{\tau(s)}=sa_{\tau(s)}^*\not= 0$. Thus $|I(\tau)|\ls 1\ls |I(\sigma)|$.
Combining this with the choice of $\sigma$, we have proved Claim 1.
\medskip

For $\pi\in S_{n-1}$ with $|I(\pi)|=1$, by $i_{\pi}$ we denote the unique element of $I(\pi)$. Without loss of generality, below we assume that
\begin{equation}\label{2.1}i_{\sigma}=\min\{i_{\pi}:\ \pi\in S_{n-1}\ \mbox{and}\ |I(\pi)|=1\}.\end{equation}
For simplicity, now we just write $i$ for $i_{\sigma}$. As $o(a_{\sigma(i)})$ divides both $i$ and $n=\exp(G)$, we have $o(a_{\sigma(i)})\mid i_n$, where
$i_n=(i,n)$.

\medskip
{\it Claim} 2: $i\mid n$.

Suppose that $i\nmid n$. Then $i_n\not=i$, $i_n\not\in I(\sigma)$ and hence $0\ne i_n a_{\sigma(i_n)}$.
Thus $o(a_{\sigma(i_n)})\nmid i_n$ and hence $o(a_{\sigma(i_n)})\nmid i$.
Therefore
$$ia_{\sigma*(ii_n)(i)}=ia_{\sigma(i_n)}\not= 0\ \ \t{and}\ \ i_n a_{\sigma*(ii_n)(i_n)}=i_n a_{\sigma(i)}= 0,$$
where $*$ is the multiplication in $S_{n-1}$ and thus $\sigma*(ii_n)$ is the product of $\sigma$ and the cyclic permutation $(ii_n)$.
So we get $|I(\sigma*(ii_n))|=1$ and $i_{\sigma*(ii_n)}=i_n<i=i_{\sigma}$, which contradicts (\ref{2.1}). This proves Claim 2.
\medskip

{\it Claim} 3: If $1\ls j<n$ and $o(a_{\sigma(j)})\nmid i$, then $i<j$ and $i\mid j$.

Assume that $1\ls j<n$ and $o(a_{\sigma(j)})\nmid i$. Then $j\not=i$ since $ o(a_{\sigma(i)})\mid i$.
For any $s=1,\ldots,n-1$ with $s\not=i,j$, we have
$$sa_{\sigma*(ij)(s)}=sa_{\sigma(s)}\not=0.$$
Also, $ia_{\sigma*(ij)(i)}=ia_{\sigma(j)}\not=0$ since $o(a_{\sigma(j)})\nmid i$.
As $|I(\sigma*(ij))|\gs |I(\sigma)|=1$, we must have $0=ja_{\sigma*(ij)(j)}=ja_{\sigma(i)}$, i.e.,  $o(a_{\sigma(i)})\mid j$.
Since $I(\sigma*(ij))=\{j\}$, we have $j=i_{\sigma*(ij)}>i=i_{\sigma}$.

Suppose that $j$ is not divisible by $i$. Then $k:=(i, j)<i$ and hence $ka_{\sigma(k)}\ne 0$ as $I(\sigma)=\{i\}$. By the
last paragraph, we must have $o(a_{\sigma(k)})\mid i$ since $k\not>i$. For any
$s=1,\ldots,n-1$ with $s\not=i,j,k$, we have
$sa_{\sigma*(kij)(s)}=sa_{\sigma(s)}\not=0$. Note that
$ia_{\sigma*(kij)(i)}=ia_{\sigma(j)}\not=0$. If $0\ne
ja_{\sigma(k)}=ja_{\sigma*(kij)(j)}$,
 then we must have $I(\sigma*(kij))=\{k\}$ and hence $i_{\sigma*(kij)}=k<i=i_{\sigma}$ which leads to a contradiction.
 Therefore, $0=ja_{\sigma(k)}$, i.e., $o(a_{\sigma(k)})\mid j$.
 Since $o(a_{\sigma(k)})$ also divides $i$,  the number $o(a_{\sigma(k)})$ must divide $(i, j)=k$, which
contradicts the fact that $ka_{\sigma(k)}\ne 0$. This proves Claim 3.

\medskip

In light of Claims 2 and 3, we have $i\in D(n)$ and
\begin{align*}|\{1\ls s<n:\ o(a_s)\nmid i\}|=&|\{1\ls j<n:\ o(a_{\sigma(j)})\nmid i\}|
\\\ls& |\{i<j<n:\ i\mid j\}|=\f ni-2.
\end{align*}
Hence, for $d=n/i\in D(n)$, we have
$$\l|\l\{1\ls s<n:\ \f nd a_s\not=0\r\}\r|<d-1,$$
which contradicts our condition (\ref{1.2}).
This proves the sufficiency. \qed

\section{Proofs of Theorems \ref{Thm1.2} and \ref{Thm1.3}}

\medskip

For a real number $x$, we let $\{x\}=x-\lfloor x\rfloor$ be its fractional part.
For any real numbers $\al$ and $\beta$, we set $\al+\beta\Z=\{\al+\beta q:\ q\in\Z\}$.

We need the following result of the second author \cite[Theorem 1]{S95}.

\addtocounter{theorem}{1}
\begin{lemma}\label{lem3.1} Let $\al_1,\ldots,\al_k$ be real numbers and let $\beta_1,\ldots,\beta_k$ be positive reals.
If $A=\{\al_s+\beta_s\Z\}_{s=1}^k$ covers consecutive
$$\bigg|\bigg\{\bigg\{\sum_{s\in I}\f1{\beta_s}\bg\}:\ I\se\{1,\ldots,k\}\bg\}\bg|$$
integers, then it covers all the integers.
\end{lemma}

\medskip
\noindent{\it Proof of Theorem \ref{Thm1.2}}. Without loss of generality, we may simply assume that $m_1,\ldots,m_{n-1}\in\{1,\ldots,n\}$.
Because Conjecture \ref{conj1.1} follows from Theorem \ref{thm1.1}, for some $\sigma\in S_{n-1}$ we have $n\nmid sm_{\sigma(s)}$ for all $s=1,\ldots,n-1$.
Note that $A=\{s+(n/m_{\sigma(s)})\Z\}_{s=1}^{n-1}$ covers $1,\ldots,n-1$ but it does not cover $0$.
By Lemma \ref{lem3.1}, the fractional parts
$$\l\{\sum_{s\in I}\f1{n/m_{\sigma(s)}}\r\}\ \ (I\se\{1,\ldots,n-1\})$$
must have more than $n-1$ distinct values. Thus, the set
$$ \l\{\sum_{i\in I}m_{i}:\ I\se\{1,\ldots,n-1\}\r\}=\l\{\sum_{s\in I}m_{\sigma(s)}:\ I\se\{1,\ldots,n-1\}\r\}$$
contains a complete system of residues modulo $n$. This concludes our proof of Theorem \ref{Thm1.2}. \qed
\medskip

To prove Theorem \ref{Thm1.3}, we need the following lemma.

\addtocounter{theorem}{1}
\begin{lemma}\label{lem3.2} {\rm (Alon's Combinatorial Nullstellensatz \cite{A})} Let $A_1,\ldots,A_n$ be finite subsets
of a field $F$ with $|A_i|>k_i$ for $i=1,\ldots,n$
where $k_1,\ldots,k_n$ are nonnegative integers.
 If the coefficient
of the monomial $x_1^{k_1}\cdots x_n^{k_n}$ in $P(x_1,\ldots,x_n)\in F[x_1,\ldots,x_n]$
is nonzero and $k_1+\cdots+k_n$ is
the total degree of $P$,
then there are $a_1\in A_1,\ldots,a_n\in A_n$ such that
$P(a_1,\ldots,a_n)\not=0$.
\end{lemma}

\medskip
\noindent{\it Proof of Theorem \ref{Thm1.3}}. Let $p$ be the smallest prime not dividing $n$. By Euler's theorem, $p^{\varphi(n)}\eq1\pmod n$, where $\varphi$ denotes
Euler's totient function. Let us consider the finite field $\mathbb F_q$ with $q=p^{\varphi(n)}$.
As $\mathbb F_q^*=\mathbb F_q\sm\{0\}$ is a cyclic group of order $q-1$, and $n$ is a divisor of $q-1$, there is an element $g\in \mathbb F_q^*$ of order $n$.
For $i=1,\ldots,n-1$ define
$$A_i:=\{g^k:\ 1\ls k\ls n-1\ \mbox{and}\ (g^k)^{m_i}\not=1\}.$$
Then $|A_i|=n-(m_i,n)\gs n-i$ for all $i=1,\ldots,n-1$.  For the polynomial
$$P(x_1,\ldots,x_{n-1}):=\prod_{1\ls i<j\ls n-1}\l(g^{a_i}x_i-g^{a_j}x_j\r),$$
we clearly have
\begin{align*}P(x_1,\ldots,x_{n-1})=&\det\l|(g^{a_i}x_i)^{j-1}\r|_{1\ls i,j\ls n-1}
\\=&\sum_{\sigma\in S_{n-1}}\mbox{sign}(\sigma)\prod_{i=1}^{n-1}\l(g^{a_i}x_i\r)^{\sigma(i)-1},
\end{align*}
where sign($\sigma$), the sign of $\sigma$, takes $1$ or $-1$ according as the permutation $\sigma$ is even or odd.
Choose $\sigma_0\in S_{n-1}$ with $\sigma_0(i)=n-i$ for all $i=1,\ldots,n-1$. Then
the coefficient of the monomial $\prod_{i=1}^{n-1}x_i^{n-1-i}$ in $P(x_1,\ldots,x_{n-1})$
coincides with
$$\mbox{sign}(\sigma_0)\prod_{i=1}^{n-1}(g^{a_i})^{n-i-1}\not=0,$$
and $\deg P=\bi {n-1}2=\sum_{i=1}^{n-1}(n-1-i)$. In view of Lemma \ref{lem3.2}, there are $x_1\in A_1,\ldots,x_{n-1}\in A_{n-1}$
such that $P(x_1,\ldots,x_{n-1})\not=0$.

Write $x_i=g^{f(i)}$ for all $i=1,\ldots,n-1$, where $f(i)\in\{1,\ldots,n-1\}$.
If $1\ls i<j\ls n-1$, then $g^{a_i+f(i)}=g^{a_i}x_i\not=g^{a_j}x_j=g^{a_j+f(j)}$ and hence
$$f(i)+a_i\not\eq f(j)+a_j\pmod n.$$
For each $i=1,\ldots,n-1$, as $(g^{f(i)})^{m_i}\not=1$ we have $n\nmid f(i)m_i$.
This completes the proof of Theorem \ref{Thm1.3}. \qed
\medskip

\Ack. The authors would like to thank the referee for helpful comments.

\end{document}